# A NUMERICAL FUNCTION IN THE CONGRUENCE THEORY

Florentin Smarandache

University of New Mexico

200 College Road

Gallup, NM 87301, USA

**Abstract**. In this paper we define a function L which will allow us to (separately or simultaneously) generalize many theorems from Number Theory obtained by Wilson, Fermat, Euler, Gauss, Lagrange, Leibniz, Moser, and Sierpinski.

**Introduction**.

1. Let A be the set $\{m \in Z / m = \pm p^{\beta}, \pm 2p^{\beta}$ with p an odd prime, $\beta \in N^*$, or $m = \pm 2^{\alpha}$ with $\alpha = 0, 1, 2$, or $m = 0\}$.

Let $m = \varepsilon p_1^{\alpha_1} \ldots p_r^{\alpha_r}$, with $\varepsilon = \pm 1$, all $\alpha_i \in N^*$, and $p_1, \ldots, p_r$ are distinct positive primes.

We construct the function L: Z x Z,

$$L(x, m) = (x + c_1) \cdot \ldots \cdot (x + c_{\varphi(m)})$$

where $c_1, \ldots, c_{\varphi(m)}$ are all modulo m rests relatively prime to m, and $\varphi$ is Euler's function.

If all distinct primes which divide x and m simultaneously are $p_{i_1}, \ldots, p_{i_r}$ then:

$$L(x, m) \equiv \pm 1 \pmod{p_{i_1}^{\alpha_{i_1}} \ldots p_{i_r}^{\alpha_{i_r}}}, \text{ when } m \in A$$

respectively $m \notin A$, and

$$L(x, m) \equiv 0 \pmod{m/(p_{i_1}^{\alpha_{i_1}} \ldots p_{i_r}^{\alpha_{i_r}})}.$$

For $d = p_{i_1}^{\alpha_{i_1}} \ldots p_{i_r}^{\alpha_{i_r}}$ and $m' = m/d$ we find

$$L(x, m) \equiv \pm 1 + k_1^0 d \equiv k_2^0 m' \pmod{m},$$

where $k_1^0, k_2^0$ constitute a particular integer solution of the diophantine equation $k_2 m' - k_1 d = \pm 1$ (the signs are chosen in accordance with the affiliation of m to A).
This result generalizes Gauss' theorem. ($c_1 \ldots c_{\varphi(m)} \equiv \pm 1 \pmod{m}$ when $m \in A$ respectively $m \notin A$) (see [1]) which generalized in its turn the Wilson's theorem (if p is prime then $(p - 1)! \equiv -1 \pmod{m}$).

<u>Proof.</u>

The following two lemmas are trivial:

**Lemma 1**. If $c_1, \ldots, c_{\varphi(p^\alpha)}$ are all modulo $p^\alpha$ rests, relatively prime to $p^\alpha$, with p an integer and $\alpha \in \mathbb{N}^*$, then for $k \in \mathbb{Z}$ and $\beta \in \mathbb{N}^*$ we have also that $kp^\beta + c_1, \ldots, kp^\beta + c_{\varphi(p^\alpha)}$ constitute all modulo $p^\alpha$ rests relatively prime to $p^\alpha$.

It is sufficient to prove that for $1 \leq i \leq \varphi(p^\alpha)$ we have $kp^\beta + c_i$ relatively prime to $p^\alpha$, but this is obvious.

**Lemma 2**. If $c_1, \ldots, c_{\varphi(m)}$ are all modulo m rests relatively prime to m, $p_i^{\alpha_i}$ divides m and $p_i^{\alpha_i+1}$ does not divide m, then $c_1, \ldots, c_{\varphi(m)}$ constitute $\varphi(m/p_i^{\alpha_i})$ systems of all modulo $p_i^{\alpha_i}$ rests relatively prime to $p_i^{\alpha_i}$.

**Lemma 3**. If $c_1, \ldots, c_{\varphi(q)}$ are all modulo q rests relatively prime to b and $(b, q) \sim 1$ then $b + c_1, \ldots, b + c_{\varphi(q)}$ contain a representative of the class $\hat{0}$ modulo q.

Of course, because $(b, q-b) \sim 1$ there will be a $c_{i_0} = q - b$, whence $b + c_{i_0} = Mq$ (multiple of q).

From this we have:

**Theorem 1**. If $(x, m/(p_{i_1}^{\alpha_{i_1}} \ldots p_{i_r}^{\alpha_{i_r}})) \sim 1$ then

$$(x + c_1) \cdot \ldots \cdot (x + c_{n(m)}) \equiv 0 \pmod{m/(p_{i_1}^{\alpha_{i_1}} \ldots p_{i_r}^{\alpha_{i_r}})}.$$

**Lemma 4**. Because $c_1 \ldots c_{\varphi(m)} \equiv \pm 1 \pmod{m}$ it results that $c_1 \ldots c_{\varphi(m)} \equiv \pm 1 \pmod{p_i^{\alpha_i}}$, for all i, when $m \in A$ respectively $m \notin A$.

**Lemma 5**. If $p_i$ divides x and m simultaneously, then $(x + c_1) \ldots (x + c_{\varphi(m)}) \equiv \pm 1 \pmod{p_i^{\alpha_i}}$, when $m \in A$ respectively $m \notin A$. Of course, from the lemmas 2 and 1, respectively 4, we have $(x + c_1) \ldots (x + c_{\varphi(m)}) \equiv$

$$\equiv c_1 \ldots c_{\varphi(m)} \equiv \pm 1 \pmod{p_i^{\alpha_i}}.$$

From the lemma 5 we obtain:

**Theorem 2**. If $p_{i_1}, \ldots, p_{i_r}$ are all primes which divide x and m simultaneously then $(x + c_1) \ldots (x + c_{\varphi(m)})$
$\equiv \pm 1 \pmod{p_{i_1}^{\alpha_{i_1}} \ldots p_{i_r}^{\alpha_{i_r}}}$, when $m \in A$ respectively $m \notin A$.

From the theorems 1 and 2 it results $L(x, m) = \pm 1 +$
$+ k_1 d = k_2 m'$, where $k_1, k_2 \in Z$. Because $(d, m') \sim 1$ the diophantine equation $k_2 m' - k_1 d = \pm 1$ admits integer solutions (the unknowns being $k_1$ and $k_2$). Hence $k_1 = m't +$
$+ k_1^0$ and $k_2 = dt + k_2^0$, with $t \in Z$, and $k_1^0$, $k_2^0$ constitute a

particular integer solution of our equation. Thus:

$$L(x, m) \equiv \pm 1 + m'dt + k_1^0 d \equiv \pm 1 + k_1^0 \pmod{m}$$

or

$$L(x, m) \equiv k_2^0 m' \pmod{m}.$$

## 2. APPLICATIONS.

(1) Lagrange extended Wilson as follows:

"if $p$ is prime, then $x^{p-1} - 1 \equiv (x + 1)(x + 2) \ldots (x + p - 1) \pmod{p}$"; we shall extend this result in the following way: For any $m \neq 0, \pm 4$ we have for $x^2 + s^2 \neq 0$ that $x^{\varphi(m_s)+s} - x^s \equiv (x + 1)(x + 2) \ldots (x + |m| - 1) \pmod{m}$, where $m_s$ and $s$ are obtained from the algorithm:

$$(0) \begin{cases} x = x_0 d_0; \quad (x_0, m_0) \sim 1 \\ m = m_0 d_0; \quad d_0 \neq 1 \end{cases}$$

$$(1) \begin{cases} d_0 = d_0^1 d_1; \quad (d_0^1, m_1) \sim 1 \\ m_0 = m_1 d_1; \quad d_1 \neq 1 \end{cases}$$

. . . . . . . . . . . . . . . . . . . . . .

$$(s-1) \begin{cases} d_{s-2} = d_{s-2}^1 d_{s-1}; \quad (d_{s-2}^1, m_{s-1}) \sim 1 \\ m_{s-2} = m_{s-1} d_{s-1}; \quad d_{s-1} \neq 1 \end{cases}$$

$$(s) \begin{cases} d_{s-1} = d_{s-1}^1 d_s; \ (d_{s-1}^1, m_s) \sim 1 \\ m_{s-1} = m_s d_s; \ d_s = 1 \end{cases}$$

(see [3] or [4]). For m a positive prime we have $m_s = m$, $s = 0$ and $\varphi(m) = m - 1$, that is Lagrange's.

(2) L. Moser enunciated the following theorem: "If p is prime, the $(p - 1)! a^p + a = Mp$", and Sierpinski (see [2], p. 57): "If p is prime then $a^p + (p - 1)! a = Mp$" which merges Wilson's and Fermat's theorems in a single one.

The function L and the algorithm from &2 will help us to generalize them too, so: if "a" and m are integers, $m \neq 0$, and $c_1, \ldots, c_{\varphi(m)}$ are all modulo m rests relatively prime to m then

$$c_1 \ldots c_{\varphi(m)} a^{\varphi(m_s)+s} - L(0, m) a^s = Mm$$

respectively

$$-L(0, m) a^{\varphi(m_s)+s} + c_1 \ldots c_{\varphi(m)} a^s = Mm,$$

or more,

$$(x + c_1) \ldots (x + c_{\varphi(m)}) a^{\varphi(m_s)+s} - L(x, m) a^s = Mm$$

respectively

$$-L(x, m) a^{\varphi(m_s)+s} + (x + c_1) \ldots (x + c_{\varphi(m)}) a^s = Mm,$$

which reunites Fermat, Euler, Wilson, Lagrange and Moser (respectively, Sierpinski).

(3) The author also obtained a partial extension of Moser's and Sierpinski's results (see [6], problem 7.140, pp. 173-174), so: if m is a positive integer, $m \neq 0, 4$, and "a" is an integer, then $(a^m - a)(m - 1)! = Mm$, reuniting Fermat and Wilson in another way.

(4) Leibniz enunciated that: "if p is prime then $(p - 2)! \equiv 1 \pmod{p}$"; we consider "$c_i < c_{i+1} \pmod{m}$" if $c_i' < c_{i+1}'$ where $0 \leq c_i' < |m|$, $0 \leq c_{i+1}' < |m|$ and $c_i \equiv c_i' \pmod{m}$, $c_{i+1} \equiv c_{i+1}' \pmod{m}$; one simply gives that if $c_1, c_2, \ldots, c_{\varphi(m)}$ are all modulo m rests relatively prime to m ($c_i < c_{i+1} \pmod{m}$ for all i, $m \neq 0$) then $c_1 c_2 \ldots c_{\varphi(m)-1} \equiv \pm 1 \pmod{m}$ if $m \in A$ respectively $m \notin A$, because $c_{\varphi(m)} \equiv -1 \pmod{m}$.